\documentclass[graybox]{Schoof_svmult}

\usepackage{type1cm}        
\usepackage{makeidx}         
\usepackage{graphicx}        
\usepackage{multicol}        
\usepackage[bottom]{footmisc}

\usepackage{newtxtext}       %
\usepackage[varvw]{newtxmath}       

\makeindex             

\usepackage{amsmath,amsfonts}%
\usepackage{amsthm}%
\usepackage{bm}
\usepackage[title]{appendix}%
\usepackage{booktabs}%
\usepackage{longtable}%
\usepackage{empheq}%
\usepackage{doi}%
\hypersetup{
  pdftitle={Title},
  allcolors=black,
  colorlinks=true
}%
\usepackage[capitalize, noabbrev]{cleveref}%
\crefname{subsection}{Subsection}{Subsections}%
\usepackage[binary-units]{siunitx}%
\DeclareSIUnit\crate{C}%
\usepackage[labelformat=simple]{subcaption}%
\usepackage[%
square,
numbers,
sort&compress,
]{natbib}
\makeatletter
\renewcommand\@biblabel[1]{#1.}
\makeatother

\usepackage{dsfont}

\def\aSi{aSi}

\makeatletter
\newcommand{\@defs@vecto}[1]{\boldsymbol{#1}}
\newcommand{\@defs@tenso}[1]{\bm{\mathbf{#1}}}

\DeclareMathOperator{\RelTol}{RelTol}
\DeclareMathOperator{\AbsTol}{AbsTol}
\DeclareMathOperator{\tr}{tr \!}

\newcommand{\de}{\,\mathrm{d}}
\newcommand{\ptl}{\partial}

\newcommand{\te}[1]{\tens{#1}}
\newcommand{\ve}[1]{\vec{#1}}

\newcommand{\trp}{\textsf{T}}

\newcommand{\grad}{\boldsymbol{\nabla}}

\newcommand{\tfinal}{t_\text{end}}

\newcommand{\halb}{\frac{1}{2}}

\newcommand{\dualreduction}{\!:\!}

\newcommand{\ton}{\text{on }}

\newcommand{\ch}{\text{ch}}
\newcommand{\el}{\text{el}}

\newcommand{\gradL}{\grad_0}
\newcommand{\divgL}{\grad_0 \!\cdot\!}
\newcommand{\something}{\boldsymbol{\cdot}}
\newcommand{\minus}{\scalebox{0.75}[1.0]{$-$}}
\newcommand{\minusone}{{\minus 1}}

\newcommand{\normalL}{\ve{n}_0}
\newcommand{\DeltaL}{\Delta_0}

\newcommand{\ext}{\text{ext}}

\newcommand{\conb}{\OB{c}}

\newcommand{\OB}[1]{\mkern 1.5mu\overline{\mkern-1.5mu#1\mkern-1.5mu}\mkern
  1.5mu}

\newcommand{\lJump}{[\![}
\newcommand{\rJump}{]\!]}

\makeatother

\begin{document}

  \title*{
    Residual Based Error Estimator for Chemical-\\
    Mechanically Coupled Battery Active Particles
  }
  \titlerunning{Residual Based Error Estimator}
  \author{Raphael Schoof
    \orcidID{\href{https://orcid.org/0000-0001-6848-3844}{0000-0001-6848-3844}},
    Lennart Fl\"ur
    \orcidID{\href{https://orcid.org/0009-0002-9728-7572}{0009-0002-9728-7572}},
    Florian Tuschner
    \orcidID{\href{https://orcid.org/0009-0004-5117-5684}{0009-0004-5117-5684}}
    and
    Willy D\"orfler
    \orcidID{\href{https://orcid.org/0000-0003-1558-9236}{0000-0003-1558-9236}}
  }
  \authorrunning{R. Schoof, L. Fl\"ur, F. Tuschner,
    W. D\"orfler}
  \institute{
    Raphael Schoof, \email{raphael.schoof@kit.edu}
    \and
    Lennart Fl\"ur, \email{lennart.fluer@student.kit.edu}
    \and
    Florian Tuschner, \email{florian.tuschner@partner.kit.edu}
    \and
    Willy D\"orfler, \email{willy.doerfler@kit.edu}
    \at
    Institute for Applied and Numerical Mathematics (IANM),
    Karlsruhe Institute of Technology (KIT),
    Englerstr.~2, 76131 Karlsruhe, Germany
  }
  \maketitle

  \abstract*{
    Adaptive finite element methods are a powerful tool to obtain numerical
    simulation results in a reasonable time.
    Due to complex chemical and mechanical couplings in lithium-ion batteries,
    numerical simulations are very helpful to investigate promising new
    battery active materials such as amorphous silicon
    featuring a higher energy density than graphite.
    Based on a thermodynamically consistent continuum model with large
    deformation and chemo-mechanically coupled approach,
    we compare
    three different spatial adaptive refinement strategies:
    \textit{Kelly-},
    \textit{gradient recovery-}
    and
    \textit{residual based}
    error estimation.
    For the residual based case,
    the strong formulation of the residual
    is explicitly derived.
    With amorphous silicon as example material,
    we investigate two 3D representative host particle geometries,
    reduced with
    symmetry assumptions to a 1D unit interval and a 2D elliptical domain.
    Our numerical studies show
    that the Kelly estimator overestimates the error,
    whereas the gradient recovery estimator leads to lower
    refinement levels
    and
    a good capture of the change of the lithium flux.
    The residual based error estimator reveals a strong dependency
    on the cell error part
    which can be improved by a more suitable choice of constants to be more
    efficient.
    In a 2D domain, the concentration has a larger influence on
    the mesh distribution than the Cauchy stress.
  }

  \abstract{
    Adaptive finite element methods are a powerful tool to obtain numerical
    simulation results in a reasonable time.
    Due to complex chemical and mechanical couplings in lithium-ion batteries,
    numerical simulations are very helpful to investigate promising new
    battery active materials such as amorphous silicon
    featuring a higher energy density than graphite.
    Based on a thermodynamically consistent continuum model with large
    deformation and chemo-mechanically coupled approach,
    we compare
    three different spatial adaptive refinement strategies:
    \textit{Kelly-},
    \textit{gradient recovery-}
    and
    \textit{residual based}
    error estimation.
    For the residual based case,
    the strong formulation of the residual
    is explicitly derived.
    With amorphous silicon as example material,
    we investigate two 3D representative host particle geometries,
    reduced with
    symmetry assumptions to a 1D unit interval and a 2D elliptical domain.
    Our numerical studies show
    that the Kelly estimator overestimates the error,
    whereas the gradient recovery estimator leads to lower
    refinement levels
    and
    a good capture of the change of the lithium flux.
    The residual based error estimator reveals a strong dependency
    on the cell error part
    which can be improved by a more suitable choice of constants to be more
    efficient.
    In a 2D domain, the concentration has a larger influence on
    the mesh distribution than the Cauchy stress.
  }



  \section{Introduction}
  \label{Schoof_sec:introduction}

  Rechargeable batteries have become an integral part of our everyday
  lives~\cite[Section~1.2]{Schoof_Writer_2019lithium-ion}, especially lithium
  (Li)-ion batteries have convinced through high energy density and long life
  time and have become one of the most popular energy storage
  system~\cite{
    Schoof_Armand_2020lithium-ion,
    Schoof_Xu_2023high-energy,
    Schoof_Zhan_2021promises,
    Schoof_Marti-Florences_2023modelling,
    Schoof_Nzereogu_2022anode}.
  Compared to the state of the art anode material
  graphite~\cite{Schoof_Zhao_2019review}, amorphous silicon~(\aSi) has the
  advantage of a nearly tenfold theoretical capacity and is therefore a very
  promising candidate for next generation Li-ion
  batteries~\cite{Schoof_Zuo_2017silicon,
    Schoof_Tomaszewska_2019lithium-ion,
    Schoof_Li_2021diverting,
    Schoof_Tian_2015high}.
  However, this benefit is accompanied by a volume change up to 300\% compared
  to
  a volume change of approximately~10\% in graphite
  anodes during lithium intake and
  extraction~\cite{Schoof_Zhang_2011review, Schoof_Zhao_2019review}. To
  investigate
  the mechanical stresses resulting from the large volume changes, efficient
  numerical simulations are necessary to understand deterioration of the
  anode, the aging process and its implications on the battery's
  lifetime,~\cite{Schoof_Miranda_2016computer,
    Schoof_Xu_2016electrochemomechanics,
    Schoof_Zhang_2011review,
    Schoof_Zhao_2019review,Schoof_Abu_2023state}
  and~\cite[Section~2.7]{Schoof_Korthauer_2018lithium-ion}.

  In order to get numerical simulation results in reasonable time, adaptive
  refinement strategies in space and time are essential, e.g.\ for phase
  separation materials like lithium iron
  phosphate~$\text{Li}_x\text{FePO}_4$~(LFP).
  This case is discussed
  in~\cite{Schoof_Castelli_2021efficient},
  where a gradient recovery estimator is used for the spatial refinement
  strategy as wells as a fully variable
  order -- variable time step size for the time integration scheme. Since there
  is no moving
  phase front in \aSi\ after the first half cycle of
  lithiation~\cite{Schoof_Di-Leo_2014cahn-hilliard-type}, no gradient motion of
  Li-concentration
  regarding the spatial coordinates can be observed as in the case of LFP.
  This makes it
  necessary to investigate whether this estimator is more efficient or not.

  In this paper, we will compare three different types of spatial refinement
  strategies. Firstly, the \textit{Kelly} error estimator, \cite{
    Schoof_Kelly_1983posteriori,
    Schoof_Arndt_2023deal}
  and~\cite[Section~4.2]{Schoof_Ainsworth_2000posteriori}, approximates the
  error
  of the jump of the solution gradients along the faces of all cells.
  This estimator is ready for use as part
  of the open source finite element
  bibliography~\texttt{deal.II}~\cite{Schoof_Arndt_2023deal}.
  Since this bibliography is used as basis for
  our implementation, we will include this estimator in our comparison.
  Secondly, the~\textit{gradient recovery} error
  estimator~\cite[Chapter~4]{Schoof_Ainsworth_2000posteriori} is used as the
  spatial
  refinement criterion as in~\cite{Schoof_Castelli_2021efficient}. Thirdly, a
  \textit{residual based} error estimator,
  see for more information \cite[Section~9.2]{
    Schoof_Brenner_2008mathematical}
  and~\cite{Schoof_Babuska_2010residual-based},
  is investigated,
  which leads to a reduction
  in the spatial error
  of the strong formulation of the respective differential equation
  to solve.
  The expressions for the residual based error indicator
  contains cell and face errors,
  whereby the inner face errors correspond in a similar way to the Kelly error
  estimator up to a chosen constant.
  See~\cite{Schoof_Ainsworth_1997posteriori,
    Schoof_Arndt_2023deal,
    Schoof_Kelly_1983posteriori}
  for the choice of the constant
  and also for further
  information~\cite[Section~1.2]{Schoof_Verfurth_1996review},
  e.g.\
  compare~\cite{Schoof_Carstensen_2003posteriori} for the $p$-laplace problem.

  The aim of this paper is the numerical comparison of these three refinement
  strategies. As part of the cell errors of the residual error estimator,
  the strong form of the underlying differential system has to be derived in
  explicit form.
  Due to the strong coupling of
  chemical and mechanical effects, the
  derivation of specific parts are presented in detail. A numerical study in 1D
  shows the differences and similarities between the three
  refinement strategies,
  whereas a 2D example considers only the gradient recovery estimator and an
  updated version of the residual based error estimator.

  The rest of this work is organized as follows: the theory of our used model
  approach is introduced in~\cref{Schoof_sec:theory} as well as the derivation
  of
  the terms of the strong formulation needed for the residual based error
  estimation. In~\cref{Schoof_sec:numerics}, the numerical solution algorithm is
  briefly presented, followed by the numerical results and discussions
  in~\cref{Schoof_sec:results}. The last~\cref{Schoof_sec:conclusion} summarizes
  our main findings and provides an outlook.



  \section{Theory}
  \label{Schoof_sec:theory}

  Firstly, this section reviews and summarizes the theory
  from~\cite{Schoof_Schoof_2022parallelization, Schoof_Schoof_2023simulation}
  for
  the formulation of our chemo-mechanically coupled model for battery active
  particles. Our theory is based on the thermodynamically consistent
  approach~\cite{Schoof_Kolzenberg_2022chemo-mechanical}. In the end, the
  derivation of all needed terms for the residual based error estimator is
  given.

\begin{figure}[!bt]
    \centering
    \includegraphics[width = 0.9\textwidth,
    page=1]{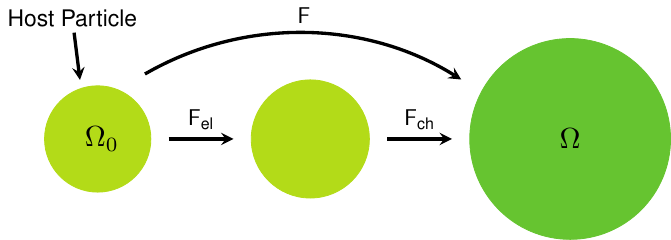}
    \caption{
      Multiplicative decomposition of the total deformation gradient~$\te{F}$
      into an elastic part~$\te{F}_\el$
      and a chemical
      part~$\te{F}_\ch$,
      compare~\cite[Figure~1]{Schoof_Schoof_2022parallelization}.}
    \label{Schoof_Figure1:deformation_theory}
  \end{figure}

  \textbf{Finite Deformation.}
  Due to the large volume changes of \aSi, we introduce the total deformation
  gradient tensor~$\te{F} \in \mathbb{R}^{d,d}$ being the derivative of the
  mapping
  from
  the referential Lagrangian domain~$\Omega_0 \subset \mathbb{R}^{d}$ to the
  current
  Eulerian domain$~\Omega \subset \mathbb{R}^{d}$ with the dimension $d \in
  \{1, 2,
  3\}$, see~\cref{Schoof_Figure1:deformation_theory} (for more information
  about the mapping,
  see~\cite[Section~2]{Schoof_Holzapfel_2010nonlinear},
  \cite[Chapter~VI]{Schoof_Braess_2007finite}
  and~\cite{Schoof_Schoof_2022parallelization, Schoof_Schoof_2023simulation,
    Schoof_Kolzenberg_2022chemo-mechanical}).
  In our case, the total deformation gradient~$\te{F}$ is completely reversible
  and can be multiplicatively split up into an elastic part~$\te{F}_\el
  (\conb, \gradL \ve{u})  \in
  \mathbb{R}^{d,d}$, caused by mechanical stress, and a chemical
  part~$\te{F}_\ch(\conb)
  \in \mathbb{R}_{\text{sym}}^{d,d}$, caused by changes in the lithium
  concentration:
  \begin{align*}
    \te{Id} + \grad_0 \ve{u}
    =
    \te{F}(\gradL \ve{u})
    =
    \te{F}_\ch(\conb)
    \te{F}_\el (\conb, \grad_0 \ve{u}).
  \end{align*}
  Both, the normalized concentration~$\conb(t, \ve{X}_0) = c(t, \ve{X}_0) /
  c_{\max} \in [0,1]$ as well as the displacement~$\ve{u}(t, \ve{X}_0) \in
  \mathbb{R}^{d}$ depend on the position vector~$\ve{X}_0 \in
  \Omega_0$ in space and the time~$t \in [0, t_\text{end}]$.
  The chemical part of the deformation gradient is defined as~$\te{F}_\ch(\conb)
  = \lambda_\ch(\conb) \te{Id}$ with $\lambda_\ch(\conb) =
  \sqrt[3]{1+v_\text{pmv} c_{\max} \conb}$
  and the elastic part~$ \te{F}_\el(\conb, \gradL \ve{u}) = \lambda_\ch^{\minus
    1}(\conb) \te{F}(\gradL \ve{u})$.
  All values for the material parameters, such as the partial molar
  volume~$v_\text{pmv}$, can be found in~\cref{Schoof_Table3:parameters}
  in~\cref{Schoof_app:simulation_parameters}.

  \textbf{Free Energy.}
  Guaranteeing a strictly positive entropy
  production~\cite{Schoof_Latz_2015multiscale, Schoof_Latz_2011thermodynamic,
    Schoof_Kolzenberg_2022chemo-mechanical, Schoof_Schammer_2021theory}
  for a thermodynamically consistent model, we base our model approach on a
  positive free energy~$\psi (\conb, \gradL \ve{u})$, considering chemical and
  mechanical effects:
  \begin{align*}
    \psi(\conb, \grad_0 \ve{u})
    = \psi_\ch(\conb) + \psi_\el(\conb, \grad_0 \ve{u}).
  \end{align*}
  The chemical part of the free energy density is formulated using an
  experimentally obtained
  open-circuit voltage (OCV)
  curve~$U_\text{OCV}$
  in \si{\volt}~\cite{
    Schoof_Chan_2007high-performance,
    Schoof_Kolzenberg_2022chemo-mechanical,
    Schoof_Keil_2016calendar,
    Schoof_Latz_2015multiscale,
    Schoof_Latz_2013thermodynamic}:
  \begin{align}
    \label{Schoof_Equation1:chemical_energy_density}
    \rho \psi_\ch(\conb)
    =
    \minus c_{\max} \int_0^{\conb} \mathrm{Fa} \, U_\text{OCV}(z) \de z
  \end{align}
  with
  the density $\rho$ of aSi and
  the Faraday constant~$\mathrm{Fa}$. We base the model for the elastic part
  on the linear elastic approach (Saint Venant--Kirchhoff model) as
  in~\cite[Section~6.5]{Schoof_Holzapfel_2010nonlinear},
  \cite[Section~VI~\S{}3]{Schoof_Braess_2007finite} and
  \cite{Schoof_Castelli_2021efficient, Schoof_Kolzenberg_2022chemo-mechanical}:
  \begin{align}
    \label{Schoof_Equation2:elastic_energy_density}
    \rho \psi_\el
    = \frac{1}{2}
    \te{E}_\el\!:\!\bm{\mathds{C}} \left[ \te{E}_\el \right] \quad
    \text{and} \quad \bm{\mathds{C}} \left[ \te{E}_\el \right]
    = \lambda \tr\,(\te{E}_\el)\te{Id} + 2 G \te{E}_\el,
  \end{align}
  the constant, isotropic stiffness tensor $\bm{\mathds{C}}$ of aSi,
  first and second Lam\'{e} constants~$ \lambda = 2 G \nu
  /\left(1-2\nu\right)$ and~$ G = E / \left[2 \left(1+2\nu\right)\right]$,
  Young's modulus~$E$ and Poisson's ratio~$\nu$. The elastic strain
  tensor~$\te{E}_\el (\conb, \gradL \ve{u}) \in \mathbb{R}_\text{sym}^{d,d}$,
  which is known in literature as
  the Green--Saint Venant strain tensor or simply \textit{the} Lagrangian strain
  tensor~\cite[Section~8.1]{Schoof_Lubliner_2006plasticity}, is defined as:
  \begin{align*}
    \te{E}_\el
    &
    =
    \frac{1}{2}
    \left(
    \te{F}_\el^{\trp}
    \te{F}_\el
    - \te{Id}\right)
    =
    \frac{1}{2}
    \left(
    \lambda_\ch^{\minus 2}
    \te{F}^{\trp}
    \te{F}
    - \te{Id}
    \right)
    =
    \frac{1}{2}
    \left(
    \lambda_\ch^{\minus 2}
    \te{C}
    - \te{Id}
    \right)
  \end{align*}
  with the right Cauchy--Green tensor~$\te{C}
  \coloneqq \te{F}^{\trp} \te{F} \in
  \mathbb{R}_\text{sym}^{d,d}$~\cite[Section~2.5]{Schoof_Holzapfel_2010nonlinear}.

  \textbf{Chemistry.}
  The change in lithium concentration inside the host material in the
  Lagrangian
  domain~$\Omega_0$
  can be described by a continuity equation~\cite{
    Schoof_Kolzenberg_2022chemo-mechanical,
    Schoof_Schoof_2022parallelization,
    Schoof_Schoof_2023simulation}:
  \begin{align*}
    \ptl_t \conb = \minus \divgL \ve{N}(\conb, \gradL \ve{u})
    \qquad
    \text{in }
    (0, \tfinal) \times \Omega_0.
  \end{align*}
  The lithium
  flux~$\ve{N}$
  is given
  as~$
  \ve{N}(\conb, \gradL \ve{u})
  = \minus m (\conb, \gradL \ve{u}) \gradL \mu(\conb, \gradL
  \ve{u})\in\mathbb{R}^{d}$
  with the mobility~$ m (\conb, \gradL \ve{u}) \coloneqq D (\ptl_c \mu (\conb,
  \gradL \ve{u}))^{\minus 1} >0$, the constant diffusion~$D >0$
  and
  the chemical
  potential~$\mu(\conb, \gradL \ve{u}) \in \mathbb{R}$,
  derived as
  the derivative of the free
  energy density
  with respect to $c$~\cite{Schoof_Kolzenberg_2022chemo-mechanical}
  as
  \begin{align}
    \label{Schoof_Equation3:chemical_potential}
    \hspace{-0.1cm}
    \mu(\conb, \gradL \ve{u})
    =
    \ptl_c (\rho \psi)(\conb, \gradL \ve{u})
    &
    =
    \minus \mathrm{Fa} \, U_\text{OCV} (\conb)
    +
    \ptl_c \te{E}_\el
    (\conb, \gradL \ve{u})
    \dualreduction
    \bm{\mathds{C}}
    [
    \te{E}_\el
    (\conb, \gradL \ve{u})
    ],
  \end{align}
  where the last
  term
  is additionally given as
  \begin{align*}
    \ptl_c\te{E}_\el
    \dualreduction
    \bm{\mathds{C}} [\te{E}_\el ]
    &
    =
    \halb
    (\minus 2)
    \lambda_\ch^{\minus 3} \ptl_c \lambda_\ch
    \te{C}
    \dualreduction
    \bm{\mathds{C}} [\te{E}_\el ]
    =
    \minus
    \lambda_\ch^{\minus 3}
    \frac{v_\text{pmv}}{3} \lambda_\ch^{\minus 2}
    \te{C}
    \dualreduction
    \bm{\mathds{C}} [\te{E}_\el ]
    \\
    &
    =
    \minus
    \frac{v_\text{pmv}}{3} \lambda_\ch^{\minus 5}
    \big(
    \te{F}^\trp \te{F}
    \big)
    \dualreduction
    \bm{\mathds{C}} [\te{E}_\el]
    =
    - \frac{v_\text{pmv}}{3 \lambda_\ch^3} (\lambda_\ch^{-2} \te{F}
    \bm{\mathds{C}} [\te{E}_\el])
    \dualreduction  \te{F}
  \end{align*}
  due to symmetry of~$\bm{\mathds{C}} [\te{E}_\el]$ and calculation rules
  for~$\square
  \dualreduction \square = \tr\,(\square^\trp \square)$. The representative host
  particle is cycled with a uniform and constant external flux~$N_\text{ext}$
  with either positive or negative sign. This external flux is applied at the
  boundary of~$\Omega_0$ and measured in terms of the
  charging-rate~(\si{C}-rate), for which we
  refer to~\cite{Schoof_Castelli_2021efficient}
  and~\cite{Schoof_Deng_2015li-ion}.
  With this definition, the simulation time $t$ and the state of charge (SOC)
  regarding the maximal capacity~\cite{Schoof_Piller_2001methods} can be related
  by
  \begin{align*}
    \mathrm{SOC}
    = \frac{c_0}{c_\text{max}} + N_\text{ext} \left[\si{C}\right] \cdot t
    [\si{\hour}].
  \end{align*}

  \textbf{Mechanics.}
  The deformation of the active material is characterized by a static balance of
  linear
  momentum~\cite{
    Schoof_Kolzenberg_2022chemo-mechanical,
    Schoof_Schoof_2022parallelization,
    Schoof_Schoof_2023simulation},
  which yields the equation:
  \begin{align*}
    \ve{0}
    =
    \minus \divgL \te{P} (\conb, \gradL \ve{u})
    \qquad
    \text{in }
    (0, \tfinal) \times \Omega_0
  \end{align*}
  with the first Piola--Kirchhoff stress tensor thermodynamical consistent
  derived as~$ \te{P}
  = \ptl_{\te{F}} (\rho \psi) = 2 \te{F} \ptl_{\te{C}} (\rho \psi) =
  \lambda_\ch^{\minus 2} \te{F} \bm{\mathds{C}} \left[ \te{E}_\el\right]\in
  \mathbb{R}^{d,d}$
  and the Cauchy
  stress~$\boldsymbol{\sigma} \in
  \mathbb{R}_{\text{sym}}^{d,d}$ in the Eulerian
  configuration$~\Omega$
  coupled
  via~$\boldsymbol{\sigma} =
  \te{P} \te{F}^{\trp} /  \det \left(\te{F}\right)$,
  compare~\cite[Sections~3.1 and~6.1]{Schoof_Holzapfel_2010nonlinear}.

  \textbf{Derivation of the strong formulation.}
  For the implementation of the residual based error estimator, both the
  divergence of the lithium flux~$\ve{N}$ and the divergence of the first
  Piola--Kirchhoff tensor~$\te{P}$ are needed in explicit form. These can be
  obtained via calculation rules
  from~\cite[Section~1.8]{Schoof_Holzapfel_2010nonlinear} for the derivatives.
  In
  a first step, $ \divgL \ve{N} (\conb, \gradL \ve{u}) $ is derived. The
  derivation of the explicit form of~$ \divgL \te{P} (\conb, \gradL \ve{u}) $
  follows subsequently.

  \textbf{Computation of }$\divgL \ve{N}$.
  With the definition of the lithium
  flux~$
  \ve{N}
  (\conb, \gradL \ve{u})
  $,
  it follows
  \begin{align}
    \divgL \ve{N}
    &
    =
    \divgL (\minus m(\conb, \gradL \ve{u}) \gradL \mu (\conb, \gradL \ve{u}))
    \nonumber
    \\
    &
    =
    \minus
    \gradL m(\conb, \gradL \ve{u}) \cdot \gradL \mu (\conb, \gradL \ve{u})
    -
    m(\conb, \gradL \ve{u}) \divgL \gradL \mu (\conb, \gradL \ve{u}))
    \nonumber
    \\
    &
    =
    \minus
    \gradL m(\conb, \gradL \ve{u}) \cdot \gradL \mu (\conb, \gradL \ve{u})
    -
    m(\conb, \gradL \ve{u}) \DeltaL \mu (\conb, \gradL \ve{u}))
    \nonumber
    \\
    &
    \hspace{0.1cm}
    =
    \minus
    \big(
    \ptl_c m(\conb, \gradL \ve{u}) \gradL c
    +
    \gradL^2 \ve{u} [\ptl_{\gradL \ve{u}} m(\conb, \gradL \ve{u})]
    \big)
    \cdot \gradL \mu (\conb, \gradL \ve{u})
    \label{Schoof_Equation4:divN}
    \\
    &
    \quad
    -
    m(\conb, \gradL \ve{u}) \DeltaL \mu (\conb, \gradL \ve{u})).
    \nonumber
  \end{align}
  \cref{Schoof_Equation4:divN} results from an application of the chain rule to
  the mobility
  with the Hessian~$\gradL^2\ve{u}(t,\ve{X}_0) \in \mathbb{R}^{d,d,d}$ of the
  displacement~$\ve{u}$, and an index reduction of the last two indices
  in~$\gradL^2\ve{u}\big[ \ptl_{\gradL \ve{u}} m(\conb, \gradL
  \ve{u})\big]\in\mathbb{R}^d$. The first term
  in~\cref{Schoof_Equation4:divN} can be formulated by using the
  definition of the mobility~$m(\conb, \gradL \ve{u}) = D (\ptl_c \mu(\conb,
  \gradL \ve{u}))^{\minus 1}$ and another application of the chain rule as
  \begin{align*}
    \ptl_c m(\conb, \gradL \ve{u})
    =
    D \ptl_c (\ptl_c \mu(\conb, \gradL \ve{u})^{\minus 1})
    =
    \minus D (\ptl_c \mu(\conb, \gradL \ve{u}))^{\minus 2}
    \ptl_c^2 \mu(\conb, \gradL \ve{u}).
  \end{align*}
  Using the definition for the chemical
  potential~$\mu(\conb, \gradL \ve{u})
  = \mu_\ch(\conb) + \mu_\el(\conb, \gradL \ve{u})$,
  compare~\cref{Schoof_Equation3:chemical_potential},
  \cref{Schoof_Equation1:chemical_energy_density}
  and
  \cref{Schoof_Equation2:elastic_energy_density},
  we have
  \begin{align*}
    \ptl_c \mu(\conb, \gradL \ve{u})
    &
    =
    \ptl_c \mu_\ch(\conb) + \ptl_c \mu_\el(\conb, \gradL \ve{u})
    =
    \ptl_c^2 (\rho \psi)_\ch(\conb)
    +
    \ptl_c^2 (\rho \psi)_\el(\conb, \gradL \ve{u})
    \\
    &
    =
    \minus \frac{1}{c_{\max}}\mathrm{Fa} \, U_\text{OCV}'
    +
    \ptl_c^2 \te{E}_\el
    \dualreduction
    \bm{\mathds{C}} [\te{E}_\el ]
    +
    \ptl_c\te{E}_\el
    \dualreduction
    \bm{\mathds{C}} [ \ptl_c \te{E}_\el ]
  \end{align*}
  and
  \begin{align*}
    \ptl_c^2 \mu(\conb, \gradL \ve{u})
    &
    =
    \ptl_c^2 \mu_\ch(\conb) + \ptl_c^2 \mu_\el(\conb, \gradL \ve{u})
    =
    \ptl_c^3 (\rho \psi)_\ch(\conb) + \ptl_c^3 (\rho \psi)_\el(\conb, \gradL
    \ve{u})
    \\
    &
    =
    \minus \frac{1}{c_{\max}^2}\mathrm{Fa} \, U_\text{OCV}''
    +
    \ptl_c^3 \te{E}_\el
    \dualreduction
    \bm{\mathds{C}} [\te{E}_\el ]
    +
    3 \ptl_c^2 \te{E}_\el
    \dualreduction
    \bm{\mathds{C}} [ \ptl_c \te{E}_\el ].
  \end{align*}
  Furthermore,
  we need
  \begin{align*}
    \ptl_{\gradL \ve{u}} m(\conb, \gradL \ve{u})
    &
    =
    D \ptl_{\gradL \ve{u}} (\ptl_c \mu(\conb, \gradL \ve{u}))^{\minusone}
    \\
    &
    =
    \minus D (\ptl_c \mu(\conb, \gradL \ve{u}))^{\minus 2}
    \ptl_{\gradL \ve{u}} \ptl_c \mu(\conb, \gradL \ve{u})
  \end{align*}
  with the last term
  \begin{align}
    \ptl_{\gradL \ve{u}} \ptl_c \mu(\conb, \gradL \ve{u})
    &
    =
    \ptl_{\gradL \ve{u}} \ptl_c \big(\mu_\ch(\conb)
    +
    \mu_\el(\conb, \gradL \ve{u})\big)
    =
    \ptl_{\gradL \ve{u}} \ptl_c \mu_\el(\conb, \gradL \ve{u})
    \nonumber
    \\
    &
    =
    \ptl_{\gradL \ve{u}}
    \big(
    \ptl_c^2 \te{E}_\el
    \dualreduction
    \bm{\mathds{C}} [\te{E}_\el ]
    +
    \ptl_c\te{E}_\el
    \dualreduction
    \bm{\mathds{C}} [ \ptl_c \te{E}_\el ]
    \big)
    \nonumber
    \\
    &
    \label{Schoof_Equation5:helping_equation1}
    =
    \ptl_{\gradL \ve{u}}
    \ptl_c^2 \te{E}_\el
    [
    \bm{\mathds{C}} [\te{E}_\el ]
    ]
    +
    \ptl_{\gradL \ve{u}} (\bm{\mathds{C}} [  \te{E}_\el]) [ \ptl_c^2\te{E}_\el ]
    \\
    &
    \label{Schoof_Equation6:helping_equation2}
    \quad
    +
    2 \ptl_{\gradL \ve{u}}
    \ptl_c
    \te{E}_\el
    [
    \bm{\mathds{C}} [ \ptl_c \te{E}_\el ]
    ].
  \end{align}
  The second term in~\cref{Schoof_Equation5:helping_equation1}
  can be computed
  with a symmetric
  tensor~$\te{S} \in \mathbb{R}_{\text{sym}}^{d,d}$
  as
  \begin{align*}
    \ptl_{\gradL \ve{u}} (\bm{\mathds{C}} [  \te{E}_\el]) [ \te{S} ]
    &
    =
    2 G \ptl_{\gradL \ve{u}} \te{E}_\el [ \te{S} ]
    +
    \lambda \te{Id} \ptl_{\gradL \ve{u}} \te{E}_\el [ \te{S} ]
    \\
    &
    =
    (2 G + \lambda)
    \ptl_{\gradL \ve{u}} \te{E}_\el [ \te{S} ]
    \\
    &
    =
    (2 G + \lambda)
    \lambda_\ch^{\minus 2} \te{F} \te{S}
  \end{align*}
  and the term
  in~\cref{Schoof_Equation6:helping_equation2}
  as
  \begin{align*}
    \ptl_{\gradL \ve{u}}
    (\ptl_c\te{E}_\el
    \dualreduction
    \bm{\mathds{C}} [ \ptl_c \te{E}_\el ])
    &
    =
    2 \ptl_{\gradL \ve{u}}\ptl_c \te{E}_\el
    [
    \bm{\mathds{C}} [\ptl_c \te{E}_\el]
    ]
    \\
    &
    =
    2 ( (\minus2) \frac{v_\text{pmv}}{3} \lambda_\ch^{\minus 5} \te{F} )\,
    (
    2 G  \ptl_c \te{E}_\el
    +
    \lambda \tr \, (\ptl_c \te{E}_\el) \te{Id}
    ).
  \end{align*}

  \textbf{Computation of }$\divgL \te{P}$.
  The divergence of the first Piola--Kirchhoff
  tensor $\te{P}(\conb, \gradL \ve{u})$
  is given as
  \begin{align*}
    \divgL \te{P}(\conb, \gradL \ve{u})
    &
    =
    \divgL
    \big(
    \lambda_\ch^{\minus 2}(\conb)
    \te{F}(\gradL \ve{u})
    \bm{\mathds{C}} [ \te{E}_\el (\conb, \gradL \ve{u})]
    \big)
    \\
    &
    =
    \te{F}
    \bm{\mathds{C}}
    [ \te{E}_\el
    ]
    \gradL \lambda_\ch^{\minus 2}
    +
    \lambda_\ch^{\minus 2}
    \divgL
    (
    \te{F}
    \bm{\mathds{C}}
    [ \te{E}_\el
    ]
    )
    \\
    &
    =
    (\minus2) \frac{v_\text{pmv}}{3} \lambda_\ch^{\minus 5}
    \te{F}
    \bm{\mathds{C}}
    [ \te{E}_\el
    ]
    \gradL c
    +
    \lambda_\ch^{\minus 2}
    \divgL
    (
    \te{F}
    \bm{\mathds{C}}
    [ \te{E}_\el
    ]
    )
    \\
    &
    =
    \ptl_{\gradL \ve{u}}\ptl_c \te{E}_\el
    [
    \bm{\mathds{C}}
    [ \te{E}_\el
    ]
    ]
    \gradL c
    +
    \lambda_\ch^{\minus 2}
    \divgL
    (
    \te{F}
    \bm{\mathds{C}}
    [ \te{E}_\el
    ]
    ).
  \end{align*}
  The final term computes to
  \begin{align*}
    \divgL
    \big(
    \te{F}(\gradL \ve{u})
    \bm{\mathds{C}} [ \te{E}_\el (c, \gradL \ve{u})]
    \big)
    &
    =
    \gradL\te{F}(\gradL \ve{u})
    [
    \bm{\mathds{C}} [ \te{E}_\el (c, \gradL \ve{u})]
    ]
    \\
    &
    \quad
    +
    \te{F}(\gradL \ve{u})
    \divgL
    (
    \bm{\mathds{C}} [ \te{E}_\el (c, \gradL \ve{u})]
    )
  \end{align*}
  with~$\gradL\te{F}(\gradL \ve{u})[\te{S}] = \gradL^2 \ve{u} [\te{S}]$ with a
  symmetric tensor~$\te{S} \in \mathbb{R}_{\text{sym}}^{d,d}$. The last term is
  turned into an explicit form
  using~$\divgL \te{A} = \tr \, (\gradL \te{A}) = \gradL \te{A} [\te{Id}]$:
  \begin{align*}
    \hspace{-0.2cm}
    \divgL
    (
    \bm{\mathds{C}}
    [ \te{E}_\el
    ]
    )
    &
    =
    2 G \divgL \te{E}_\el (\conb, \gradL \ve{u})
    +
    \lambda \divgL
    \big(
    \tr \, (\te{E}_\el (\conb, \gradL \ve{u})) \te{Id}
    \big)
    \\
    &
    =
    2 G \divgL
    \Big(
    \halb (\lambda_\ch^{\minus 2}(\conb) \te{C}(\gradL \ve{u}) -
    \te{Id})
    \Big)
    \\
    &
    \quad
    +
    \lambda \divgL
    \Big(
    \tr \, (\halb (\lambda_\ch^{\minus 2}(\conb) \te{C}(\gradL \ve{u}) -
    \te{Id}))
    \te{Id}
    \Big)
    \\
    &
    =
    2 G \halb
    \big(
    \te{C}
    \gradL \lambda_\ch^{\minus 2}
    +
    \lambda_\ch^{\minus 2}
    \divgL
    \te{C}
    \big)
    +
    \lambda \halb
    \divgL
    \big(
    \lambda_\ch^{\minus 2}
    \tr \, ( \te{C}
    )
    \te{Id}
    \big)
    \\
    &
    =
    G
    \big(
    \te{C}
    (\minus 2)\lambda_\ch^{\minus 3} \ptl_c \lambda_\ch \gradL c
    +
    \lambda_\ch^{\minus 2}
    \gradL \te{C}
    [\te{Id}]
    \big)
    \\
    &
    \quad
    +
    \lambda \halb
    \te{Id}
    \gradL
    \big(
    \lambda_\ch^{\minus 2}
    \tr \, ( \te{C}
    )
    \big)
    + \ve{0}
    \\
    &
    =
    G
    \big(
    \te{C}
    (\minus 2)\lambda_\ch^{\minus 3}
    \frac{v_\text{pmv}}{3}\lambda_\ch^{\minus 2} \gradL c
    +
    \lambda_\ch^{\minus 2}
    \gradL^2 \ve{u} [\ptl_{\grad \ve{u}} \te{F} [\ptl_{\te{F}} \te{C}
    [\te{Id}]]]
    \big)
    \\
    &
    \quad
    +
    \lambda \halb
    \big(
    \tr \, ( \te{C}
    )
    \gradL
    \lambda_\ch^{\minus 2}
    +
    \lambda_\ch^{\minus 2}
    \gradL \tr \, ( \te{C}
    )
    \big)
    \\
    &
    =
    G
    \big(
    2
    \ptl_c \te{E}_\el
    \gradL c
    +
    \lambda_\ch^{\minus 2}
    \gradL^2 \ve{u} [ 2\te{F} ]
    \big)
    \\
    &
    \quad
    +
    \halb
    \lambda
    \big(
    \tr \, ( \te{C}
    )
    (\minus 2)
    \lambda_\ch^{\minus 3} \frac{v_\text{pmv}}{3} \lambda_\ch^{\minus 2} \gradL
    c
    +
    \lambda_\ch^{\minus 2}
    \gradL^2 \ve{u}
    [\ptl_{\grad \ve{u}} \te{F} [\ptl_{\te{F}} \te{C} [\te{Id}]]]
    \big)
    \\
    &
    =
    2 G
    \big(
    \ptl_c \te{E}_\el
    \gradL c
    +
    \lambda_\ch^{\minus 2}
    \gradL^2 \ve{u} [ \te{F} ]
    \big)
    +
    \lambda
    \big(
    \tr \, ( \ptl_c \te{E}_\el
    )
    \gradL c
    +
    \lambda_\ch^{\minus 2}
    \gradL^2 \ve{u}
    [\te{F}]
    \big)
    \\
    &
    =
    \bm{\mathds{C}}
    [
    \ptl_c \te{E}_\el
    ]
    \gradL c
    +
    (
    2G + \lambda
    )
    \lambda_\ch^{\minus 2}
    \gradL^2 \ve{u} [ \te{F} ].
  \end{align*}
  Now, we have all terms to apply the residual based error estimator.



  \section{Numerical Solution Procedure}
  \label{Schoof_sec:numerics}
  Now, all important aspects for the numerical treatment
  and the three adaptive refinement possibilities,
  being numerically compared,
  are shortly introduced.

  \textbf{Problem Formulation.}
  Firstly, we refer
  to~\cref{Schoof_Table2:normalization}
  for the normalization of our used model parameters.
  From now on, we consider only the dimensionless variables without any
  accentuation.
  Secondly, we obtain our problem formulation
  of the dimensionless initial boundary value problem
  from~\cref{Schoof_sec:theory}:
  let~$t_\text{end}>0$ the final simulation time
  and~$\Omega_0$ a bounded electrode particle as reference configuration with
  dimension~$d = 3$.
  Find the normalized
  concentration~$c
  \colon [0, t_\text{end}] \times \overline{\Omega}_0
  \rightarrow [0,1]$,
  the chemical
  potential~$\mu
  \colon [0, t_\text{end}] \times \overline{\Omega}_0
  \rightarrow \mathbb{R}$
  and the displacement~$\ve{u}
  \colon [0, t_\text{end}] \times \overline{\Omega}_0
  \rightarrow \mathbb{R}^{d}$ satisfying

  \begin{subequations}
    \label{Schoof_Equation7:problem_ibvp}
    \begin{empheq}[left=\empheqlbrace]{alignat=4}
      \partial_t c
      &
      =
      \minus \divgL \ve{N}
      (c, \grad_0 \ve{u})
      &&
      \qquad
      &&
      \text{in } (0, \tfinal) \times
      \Omega_0,
      \label{Schoof_Equation7:problem_ibvp_a}
      \\
      \mu
      (c, \grad_0 \ve{u})
      &
      =
      \partial_c (\rho \psi)(c, \grad_0 \ve{u})
      &&
      &&
      \text{in } (0, \tfinal)
      \times \Omega_0,
      \\
      \ve{0}
      &
      =
      \minus \divgL \te{P}(c, \grad_0\ve{u})
      &&
      &&
      \text{in } (0,
      \tfinal) \times \Omega_0,
      \\
      \ve{N}(c, \grad_0\ve{u}) \cdot \normalL
      &
      =
      N_\text{ext}
      &&
      &&
      \text{on } (0, \tfinal)
      \times \partial \Omega_0,
      \\
      \te{P}(c, \grad_0\ve{u}) \cdot \normalL
      &
      =
      \ve{0}
      &&
      &&
      \text{on } (0, \tfinal)
      \times \partial \Omega_0,
      \\
      c(0, \something)
      &
      =
      c_0
      &&
      &&
      \text{in } \Omega_0
      \label{Schoof_Equation7:problem_ibvp_f}
    \end{empheq}
  \end{subequations}
  with an initial
  condition~$c_0$
  being consistent with the boundary conditions.
  Our quantities of interest are the concentration~$c$
  and the Cauchy stress~$\sigma$, computed in a postprocessing stress with the
  help of the three solution variables.
  In the case of lithiation of the host particle,
  the sign of the external lithium
  flux~$N_\text{ext}$
  is positive.
  Contrarily, in the case of delithiation,
  the sign of the external lithium flux is negative.
  Appropriate boundary conditions for the displacement are used to exclude rigid
  body motions.
  The original formulation of the chemical deformation
  gradient~$\te{F}_\ch$
  is written for
  the three-dimensional case, however, it is applicable and mathematically valid
  also for all variables and
  equations in lower dimensions as well.

  \textbf{Numerical Solution Procedure.}
  We apply the finite element method on
  the model equations
  of~\cref{Schoof_Equation7:problem_ibvp_a}\,--\,\cref{Schoof_Equation7:problem_ibvp_f}
  for the spatial
  discretization similar to~\cite{Schoof_Schoof_2023simulation}
  and~\cite{Schoof_Castelli_2021efficient}. Furthermore, we use the family of
  numerical differentiation formulas (NDFs) in a variable-step, variable-order
  algorithm~\cite{Schoof_Reichelt_1997matlab,
    Schoof_Shampine_1997matlab,
    Schoof_Shampine_1999solving,
    Schoof_Shampine_2003solving},
  because the properties of the resulting differential algebraic
  equation~$(\text{DAE})$ are similar to stiff ordinary differential
  equations~\cite{Schoof_Castelli_2021efficient}. The change in time step
  sizes~$\tau_n$ and order~$k_n$ is handled by an error control, see Algorithm~1
  in~\cite{Schoof_Castelli_2021efficient}. The adaptivity of Algorithm~1 is then
  achieved with a mixed error control using thresholds~$\RelTol_t$, $\AbsTol_t$,
  $\RelTol_x$ and~$\AbsTol_x$, compare~\cite{Schoof_Castelli_2021efficient}
  and~\cite[Section~1.4]{Schoof_Shampine_2003solving}. Finally, the resulting
  nonlinear system is linearized with the Newton--Raphson method and is solved
  for the updates with a direct LU-decomposition. In this work, we compare three
  error estimators for the local criterion of the adaptive mesh refinement and
  coarsening: the \textit{Kelly-},
  the \textit{gradient recovery-}
  and the \textit{residual based}
  error estimator,
  as described in the introduction.
  For later use, we state the definition of the residual:
  \begin{align*}
    \ve{R}_h^{n+1}
    =
    \begin{pmatrix}
      \frac{1}{\hat{\tau}_n} \te{M} (\ve{y}^{n+1} -\ve{\Phi}^{n})
      +
      \divgL
      \ve{N}_h^{n +1}
      \\
      \mu_h^{n+1} - \partial_c (\rho\psi)_h^{n+1}
      \\
      \divgL \te{P}_h^{n+1}
    \end{pmatrix}
  \end{align*}
  with the mass matrix~$\te{M}$
  of the finite element space of
  dimension~$N$
  regarding the
  concentration~$c$,
  the discrete solution
  vector~$\ve{y} \colon [0, t_{\text{end}}] \rightarrow \mathbb{R}^{(2 + d)N}$
  with~$\ve{y}^{n+1} \approx \ve{y}(t_{n+1})$,
  containing all three discrete
  components~$c_h$,
  $\mu_h$
  and~$\ve{u}_h$.
  $\hat{\tau}_n$ is a modified time step size with coefficients of the selected
  order~$k_n$
  at time~$t_n$
  and $\ve{\Phi}^{n}$ consists of the solutions of some former time
  steps~$\ve{y}^{n}, \dots,
  \ve{y}^{n-k_n}$~\cite[Section~2.3]{Schoof_Shampine_1997matlab}.

  In total, the residual based error estimation~$\text{est}_x$ can be split up
  into a cell fraction and a face fraction,
  whereby the latter one includes also the
  boundary parts. For the first one, we sum up the~$L^2$-norm of the residual
  over all cells~$K$ of our triangulation, weighted with the largest diagonal of
  the cell~$h_K$ in~\cref{Schoof_Equation8:residual_cell_part}. For the latter
  one,
  we
  add up all jumps~$\lJump \something \rJump$ of the function in square brackets
  at the inner faces~$\mathcal{F}^{\mathrm{o}}$ of all cells
  in~\cref{Schoof_Equation9:residual_faces_inner_part} and additionally
  as well as the boundary conditions for all faces~$\mathcal{F}^{\ptl
    \Omega_{\text{0}}}$ at all boundary cells
  in~\cref{Schoof_Equation10:residual_faces_boundary_part}. Both face parts are
  weighted
  with~$h_{K}/24$ as proposed
  in~\cite[Section~4.2]{Schoof_Ainsworth_2000posteriori}.
  In summary, we have the
  overall error estimator with constants~$ \gamma_1 >0$ and~$ \gamma_2 > 0$
  on one cell~$K$
  computed as:
  \begin{align}
    \vspace{-1.0cm}
    \label{Schoof_Equation8:residual_cell_part}
    \text{est}_K^2
    &
    \coloneqq
    \gamma_1 \eta_\text{cell,K}^2 + \gamma_2 \eta_\text{face,K}^2
    \\
    &
    =
    \gamma_1 h_K^2
    \| \textbf{R}_h^{n+1} \|_{L^2(K)}^2
    \nonumber
    \\
    &
    \label{Schoof_Equation9:residual_faces_inner_part}
    \quad
    +
    \gamma_2
    \frac{h_K}{24}
    \sum_{F \in \mathcal{F}^{\mathrm{o}} \cap \ptl K}
    \Bigl(
    \|
    \lJump
    \ve{N}_h^{n+1} \cdot \normalL
    \rJump \|_{L^2(F)}^2
    +
    \|
    |
    \lJump
    \te{P}_h^{n+1} \cdot \normalL
    \rJump
    |_2^2
    \|_{L^2(F)}^2
    \Bigr)
    \\
    &
    \label{Schoof_Equation10:residual_faces_boundary_part}
    \quad
    +
    \gamma_2
    \frac{h_K}{24}
    \sum_{F \in \mathcal{F}^{\ptl \Omega_{\text{0}}} \cap \ptl K}
    \Bigl(
    \|
    \ve{N}_h^{n+1} \cdot \normalL - N_\text{ext}
    \|_{L^2(F)}^2
    +
    \|
    |
    \te{P}_h^{n+1} \cdot \normalL - \ve{0}
    |_2^2
    \|_{L^2(F)}^2
    \Bigr).
  \end{align}
  The total error~$\text{est}_x$ is then given as $\text{est}_x^2=\sum_{K}
  \text{est}_K^2$.
  Keep in mind that we have to add further boundary conditions for some
  additional
  artificial boundaries for the specified computational domain as explained
  in~\cref{Schoof_subsec:simulation_setup}.



  \section{Numerical Studies}
  \label{Schoof_sec:results}
  In this section, we specify our simulation setup,
  present
  and discuss
  our simulation results
  of the comparison
  of the three spatial refinement strategies for the presented model
  of~\cref{Schoof_sec:theory}.

  \subsection{Simulation Setup}
  \label{Schoof_subsec:simulation_setup}
  We choose \aSi\ as host material for our derived theory with material
  parameters and dimensionless values given in~\cref{Schoof_Table3:parameters}.
  The
  used experimental~$\text{OCV}$-curve is taken
  from~\cite[Equation~(SI-51)]{Schoof_Kolzenberg_2022chemo-mechanical} and
  displayed in~\cite[Figure~2]{Schoof_Schoof_2022parallelization}. For
  lithiation,
  we apply an external constant lithium flux~$N_\text{ext} = 1\si{C}$ and for
  delithiation, the value~$N_\text{ext} = \minus 1\si{C}$ is given.
  Following~\cite{
    Schoof_Kolzenberg_2022chemo-mechanical,
    Schoof_Schoof_2023simulation,
    Schoof_Schoof_2023efficient},
  we start with a constant concentration~$c_0 = 0.02$ and determine one cycle
  period of lithiation or delithiation with a duration $0.90\,\si{h}$, so we are
  in the range of $\text{SOC} \in [0.02, 0.92]$ and time~$t \in [0.0~\si{h},
  0.90\,\si{h}]$ for lithiation, followed by a delithiation with the same
  duration of~$0.90~\si{h}$ and decreasing~$\text{SOC}$.
  After the discharging, another charging cycle is added and so on.
  For a total simulation time of~$t=2.7~\si{h}$, we have following cycles:
  charging--discharging--charging.

\begin{figure}[!t]
    \centering
    \includegraphics[width = 0.75\textwidth,
    page=1]{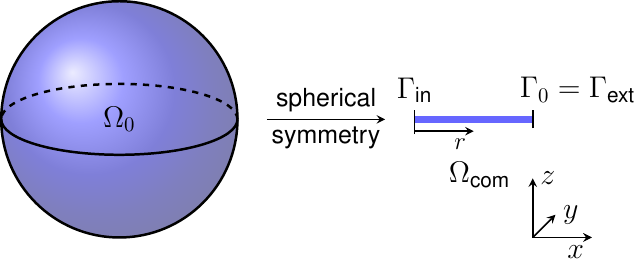}
    \caption{Dimension reduction of a 3D unit sphere
      to the 1D unit interval with spherical symmetry, based
      on~\cite[Figure~B.1]{Schoof_Castelli_2021numerical}.}
    \label{Schoof_Figure2:particle_1d}
  \end{figure}

\begin{figure}[!b]
    \centering
    \includegraphics[width = 0.99\textwidth,
    page=2]{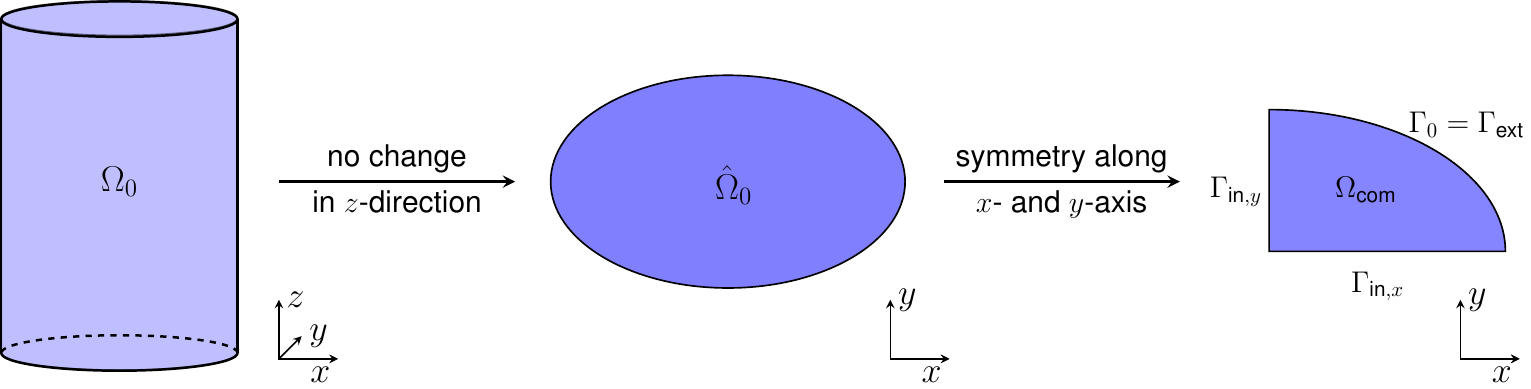}
    \caption{Dimension reduction of a 3D elliptical nanowire
      to the 2D quarter ellipsoid, based
      on~\cite[Figure~5]{Schoof_Schoof_2023simulation}.
    }
    \label{Schoof_Figure3:particle_2d}
  \end{figure}

  \textbf{Geometrical Setup.}
  Similar
  to~\cite{Schoof_Schoof_2023efficient},
  we choose a representative 3D spherical particle with domain~$\Omega_0$,
  reduced to a 1D computational domain~$\Omega_\text{com}$ as
  shown in~\cref{Schoof_Figure2:particle_1d} in terms of the radial
  variable~$r$, and a quarter ellipsoid as computational
  domain~$\Omega_\text{com}$ resulting from a 3D nanowire domain~$\Omega$ with
  no
  changes in $z$-direction and symmetry around the $x$- and $y$-axis as
  in~\cref{Schoof_Figure3:particle_2d}. For the adaption of the quadrature
  weight and
  the initial conditions, we refer to the literature cited.

  However, we have to mention the additional artificial inner boundary
  conditions,
  indicated with the subscript~$\square_\text{in}$.
  For the 1D computational domain, we have
  \begin{align*}
    \ve{N} \cdot \ve{n}_0 = 0,
    \quad
    u = 0,
    \quad
    \ton
    (0, t_\text{end} \times \Gamma_\text{in})
  \end{align*}
  and
  for the 2D computational domain
  \begin{alignat*}{3}
    \ve{N} \cdot \ve{n}_0 &= 0,
    \quad
    u_y
    &&= 0,
    \quad
    &&\ton
    (0, t_\text{end} \times \Gamma_{\text{in},x}),
    \\
    \ve{N} \cdot \ve{n}_0 &= 0,
    \quad
    u_x
    &&= 0,
    \quad
    &&\ton
    (0, t_\text{end} \times \Gamma_{\text{in},y}).
  \end{alignat*}
  These new terms also appear in the computation
  of the boundary face terms of the residual base error estimator
  of~\cref{Schoof_Equation10:residual_faces_boundary_part}.
  For the 1D computational domain,
  we have to add
  \begin{align*}
    \gamma_2
    \frac{h_K}{24}
    \sum_{F \in \mathcal{F}^{\Gamma_{\text{in}}} \cap K}
    \Bigl(
    \|
    \ve{N}_h^{n+1} \cdot \normalL - 0
    \|_{L^2(F)}^2
    +
    \|
    u^{n+1} - 0
    \|_{L^2(F)}^2
    \Bigr)
  \end{align*}
  to the already mentioned face terms
  and for the 2D computational domain
  \begin{align*}
    &
    \gamma_2
    \frac{h_K}{24}
    \sum_{F \in \mathcal{F}^{\Gamma_{\text{in},x}} \cap K}
    \Bigl(
    \|
    \ve{N}_h^{n+1} \cdot \normalL - 0
    \|_{L^2(F)}^2
    +
    \|
    u_y^{n+1} - 0
    \|_{L^2(F)}^2
    \Bigr)
    \\
    +
    &
    \gamma_2
    \frac{h_K}{24}
    \sum_{F \in \mathcal{F}^{\Gamma_{\text{in},y}} \cap K}
    \Bigl(
    \|
    \ve{N}_h^{n+1} \cdot \normalL - 0
    \|_{L^2(F)}^2
    +
    \|
    u_x^{n+1} - 0
    \|_{L^2(F)}^2
    \Bigr).
  \end{align*}

  \textbf{Implementation Details.}
  We use an isoparametric fourth-order
  Lagrangian finite element method with isoparametric
  representation of the curved boundary for our numerical simulations. The
  \texttt{C++} finite element
  library~\texttt{deal.II}~\cite{Schoof_Arndt_2023deal} is taken as basis as
  well
  as the interface to the Trilinos
  library\cite[Version~12.8.1]{Schoof_Team_2023trilinos} and the UMFPACK
  package~\cite[Version~5.7.8]{Schoof_Davis_2004algorithm} for the
  LU-decomposition. All simulations for the 1D computational domain are
  performed
  on a desktop computer with \SI{64}{\giga\byte}~RAM, Intel~i5-9500~CPU, GCC
  compiler version~10.5 and operating system Ubuntu 20.04.6 LTS, whereas the
  simulations for the 2D computational domain are executed with GCC compiler
  version~12.1 with a single node of the BwUniCluster~2.0 with 40 Intel Xeon
  Gold~6230 with~\SI{2.1}{\giga\hertz} and \SI{96}{\giga\byte}~RAM,
  on~\cite{Schoof_Hpc-team_2024bwunicluster-2-0-hardware-and-architecture}. The
  highly resolved solution, which is used to compute the error
  in~\cref{Schoof_Table1:norm_comparison}, is also computed on BwUniCluster~2.0.
  Unless otherwise stated, we set for the space and time adaptive algorithm
  following parameters: tolerances~$\RelTol_t = \RelTol_x = \num{1e-5}$,
  $\AbsTol_t = \AbsTol_x = \num{1e-8}$, initial time step size~$\tau_0 =
  \num{1e-6}$, final simulation time~$t_\text{end}=2.7$, maximal time step
  size~$\tau_{\max} = 0.1$, number of initial refinements~7, minimal refinement
  level~3 for the 1D computational domain and level~1 for the 2D computational
  domain, maximal refinement level of 20 and marking parameters for the local
  mesh coarsening and refinement, $\theta_{\mathrm{c}} = 0.05$ and
  $\theta_{\mathrm{r}}=0.3$. Furthermore, we use automatic differentiation~(AD)
  for the computation of the Newton
  matrix~\cite{Schoof_Schoof_2023efficient}
  and OpenMP Version~4.5 for shared memory parallelization for assembling
  the
  system matrices,
  residuals
  and the Kelly- and gradient recovery estimator.

  \subsection{Numerical Results}
  \label{Schoof_subsec:numerical_results}

  \begin{figure}[t]
    \centering
    \includegraphics[width = 0.99\textwidth,
    page=1]{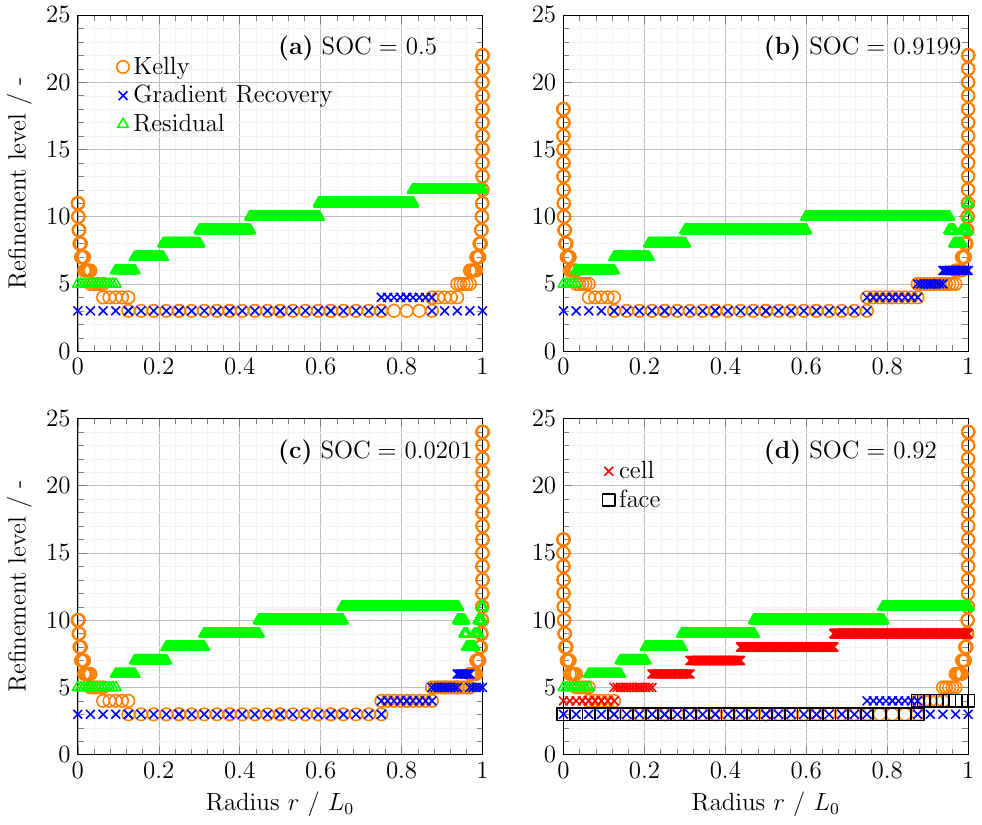}
    \caption{
      Comparison of three different refinement strategies
      (Kelly error estimator,
      gradient recovery error estimator
      and residual based error estimator
      with~$\gamma_1 = 1.0$
      and~$\gamma_2 = 1.0$)
      over the particle radius
      at different~$\text{SOC}$:
      $\text{SOC} = 0.5$
      in~(a),
      $\text{SOC} = 0.9199$
      in~(b),
      $\text{SOC} = 0.0201$
      in~(c)
      and
      $\text{SOC} = 0.92$
      in~(d),
      whereas the last one contains also the two fractions of the cell and the
      face error of the residual based error estimator.
    }
    \label{Schoof_Figure4:residual_level_1d}
  \end{figure}

  \textbf{Example 1: 1D Spherical Symmetry.}
  In a first step,
  we compare the three adaptive refinement strategies at four different time
  steps during our simulation with $t_\text{end} = 2.7~\si{h}$.
  \cref{Schoof_Figure4:residual_level_1d}
  shows the refinement levels
  for the Kelly error estimator,
  the gradient recovery estimator
  and
  the residual based error estimator
  with~$\gamma_1 = 1.0$
  and~$\gamma_2=1.0$
  over the particle radius
  at~$t=0.48~\si{h}$
  corresponding
  to~$\text{SOC} = 0.50$
  in Subfigure~(a),
  at~$t=0.9001~\si{h}$
  corresponding
  to~$\text{SOC} = 0.9199$
  in Subfigure~(b),
  at~$t=1.8001~\si{h}$
  corresponding
  to~$\text{SOC} = 0.0201$
  in Subfigure~(c)
  and the final simulation time
  at~$t_\text{end}=2.7~\si{h}$
  corresponding
  to~$\text{SOC} = 0.92$
  in Subfigure~(d).
  For the simulation with the Kelly error estimator,
  we have increased the maximal refinement level to~25 and also updated the
  tolerances
  to~$\RelTol_t = \RelTol_x = \num{1e-4}$,
  $\AbsTol_t = \AbsTol_x = \num{1e-7}$
  to get a stable simulation.
  Is is clearly evident
  that the Kelly error estimator over-refines the boundary parts,
  since the estimator has no information about the exact boundary terms.
  The gradient recovery estimator has only at some range
  around~$
  r= 0.8$
  a higher refinement level
  in~\cref{Schoof_Figure4:residual_level_1d}(a)
  and~\cref{Schoof_Figure4:residual_level_1d}(d),
  whereas shortly after the sign change of the external lithium flux also a
  higher
  refinement level is reached at the particle surface
  around~$
  r= 1.0$
  in the other two subfigures.
  The residual based error estimator features some stair-like behavior with
  increasing refinement level from smaller to larger radius values
  in~\cref{Schoof_Figure4:residual_level_1d}(a)
  and~\cref{Schoof_Figure4:residual_level_1d}(d).
  Shortly after the change of the sign of the lithium flux
  in~\cref{Schoof_Figure4:residual_level_1d}(b)
  and~\cref{Schoof_Figure4:residual_level_1d}(c),
  there is a drop of the refinement level at the particle surface.
  This could be explained by the large influence of the cell error fraction
  on the residual based refinement strategy,
  compare the two cell and face fractions,
  which are
  additionally added
  in~\cref{Schoof_Figure4:residual_level_1d}(d).
  Due to the change of the lithium flux,
  the cell fraction error is smaller
  (since the lithium flux direction close to particle surface
  around~$
  r=1.0$ has also changed)
  and therefore, a smaller refinement level is sufficient.
  The behavior of the face error of the residual based strategy
  in~\cref{Schoof_Figure4:residual_level_1d}(d)
  leads to a similar refinement as the gradient recovery strategy.

  \begin{figure}[t]
    \centering
    \includegraphics[width = 0.99\textwidth,
    page=1]{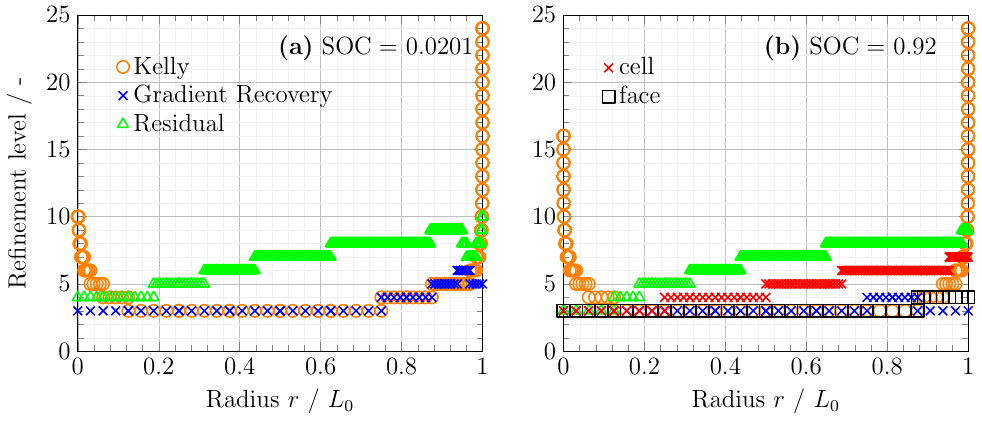}
    \caption{
      Same situation as in~\cref{Schoof_Figure4:residual_level_1d}(c)-(d), but
      with
      updated version of the residual based error estimator with
      $\gamma_1 = 0.001$ and $\gamma_2=1.0$.
    }
    \label{Schoof_Figure5:residual_level_1d_modified}
  \end{figure}

  \begin{figure}[b]
    \centering
    \includegraphics[width = 0.99\textwidth,
    page=1]{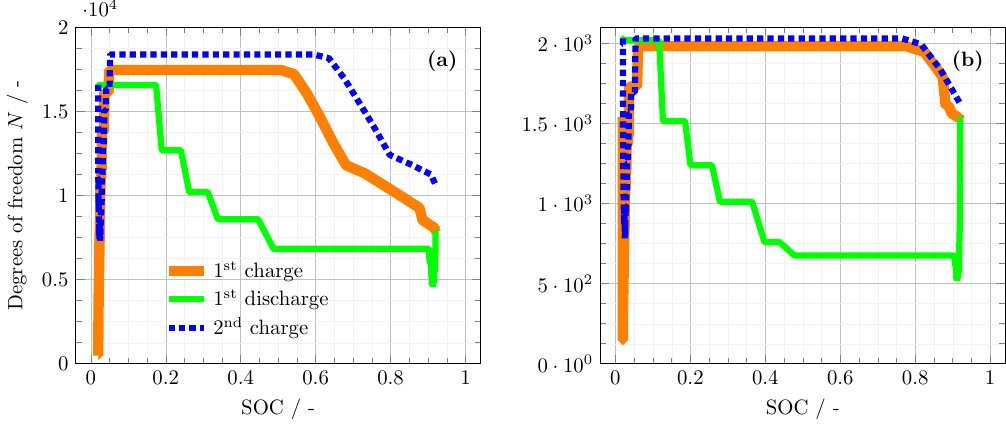}
    \caption{
      Number of DOFs
      over $\text{SOC}$
      for the residual based error estimator
      with~$\gamma_1 = 1.0$ and $\gamma_2=1.0$
      in~(a)
      and
      updated version with $\gamma_1 = 0.001$ and $\gamma_2=1.0$
      in~(b).
    }
    \label{Schoof_Figure6:residual_dofs_1d}
  \end{figure}

  It is known in literature
  that the edge fraction is the dominating one
  and the more important one for the error estimation,
  see for
  example~\cite{Schoof_Carstensen_1999edge}
  and~\cite[Section~1.6]{
    Schoof_Verfurth_2013posteriori}.
  Therefore,
  we have additionally performed some simulations with
  an updated version of the residual based error estimator
  with~$\gamma_1 = 0.001$
  and~$
  \gamma_2=1.0$.
  Our choice of~$\gamma_1$ is suitable for reducing the influence of the
  cell proportion,
  but requires further efforts,
  e.g. possible influence of the chosen finite element order
  or
  a possible comparison in the order of magnitude to the computed numerical
  error.
  However, these considerations will be postponed to future work.
  The smaller fraction of the cell error clearly reduces the refinement level of
  the residual based refinement strategy,
  comparing~\cref{Schoof_Figure5:residual_level_1d_modified}(a)-(b)
  to~\cref{Schoof_Figure4:residual_level_1d}(c)-(d).

  This influence of the cell error fraction is further investigated
  in~\cref{Schoof_Figure6:residual_dofs_1d},
  where the number of degrees of freedom (DOFs)~$(2+d)N$
  is plotted for the residual based error estimator
  with~$
  \gamma_1 = 1.0$
  and~$
  \gamma_2=1.0$.
  in Subfigure~(a)
  and the updated version
  with~$
  \gamma_1 = 0.001$
  and~$
  \gamma_2=1.0$
  in Subfigure~(b).
  Both plots show some qualitative similarities,
  which are the steep increase of the number of DOFs
  after the start of the first lithiation
  followed by some plateau and a decrease in the end of the first charge.
  The first discharge looks also similar,
  since the number of DOFs gets smaller but reaches a higher level at the end
  of the first discharge.
  For the second charge,
  there is the same behavior in both plots visible
  featuring a large drop in the number of~DOFs and following then the same way
  as
  the first charge.
  However, there is one difference:
  the absolute number of DOFs is reduced by a factor of~$10$
  in the updated version.
  This shows again the significant influence of the cell error fraction.

  \begin{table}[t]
    \centering
    \caption{
      Error between the numerical solution compared to a highly
      resolved numerical solution regarding the respective norm:
      $\|\ve{y} - \ve{y}_\text{highly\_resolved} \|_{\square}$.
    }
    \label{Schoof_Table1:norm_comparison}
    \begin{tabular}[t]{@{}lccc@{}}
      \toprule
      \textbf{Error estimator}
      &
      DOFs
      &
      $L^2\text{-Norm}$
      &
      $H^1\text{-Norm}$
      \\
      \midrule
      Gradient recovery
      &$\num{111}$
      &$\num{1.20e-06}$
      &$\num{3.02e-06}$
      \\[\defaultaddspace]

      Residual
      &$\num{15939}$
      &$\num{3.08e-07}$
      &$\num{1.07e-06}$
      \\[\defaultaddspace]

      Residual updated
      &$\num{1575}$
      &$\num{5.44e-07}$
      &$\num{1.65e-06}$
      \\[\defaultaddspace]

      \bottomrule
    \end{tabular}
  \end{table}

  In~\cref{Schoof_Table1:norm_comparison},
  we display the~$L^2\text{-Norm}$
  and~$H^1\text{-Norm}$
  of the numerical solution~$\ve{y}$
  compared with a highly resolved solution~$\ve{y}_\text{highly\_resolved}$
  at time $t=0.2~\si{h}$.
  For the numerical results of the table,
  we use the exact computation of the Newton matrix and no AD.
  The high resolved solution is computed with~$N=1.5$ million DOFs
  and
  updated~$\RelTol_t = \num{1e-11}$
  and~$\AbsTol_t = \num{1e-14}$.
  The numerical solution~$\ve{y}$ is computed
  with~$\RelTol_t = \num{1e-7}$,
  $\AbsTol_t = \num{1e-10}$,
  initial time step size~$\tau_0 = \num{1e-10}$
  and initial refinement level with 17.
  We compare the results of the gradient recovery-, the residual based- and the
  updated residual based error estimator.
  On the one hand,
  it can be seen
  that the gradient recovery displays the lowest
  number of DOFs but also the largest errors.
  On the other hand, residual based methods feature smaller norm errors but have
  also significant larger number of DOFs.
  Interestingly,
  the updated case leads to a factor of ten fewer number of DOFs,
  whereas the error norms are at a comparable level.
  From this,
  we can conclude that
  it is more efficient to use the updated version of the residual based error
  estimator due to the smaller number of DOFs.
  Further,
  we are accurate enough with the updated version
  and do not need to solve the numerical simulations for the solutions with the
  higher number of DOFs.

  \begin{figure}[t]
    \centering
    \includegraphics[width = 0.65\textwidth,
    page=1]{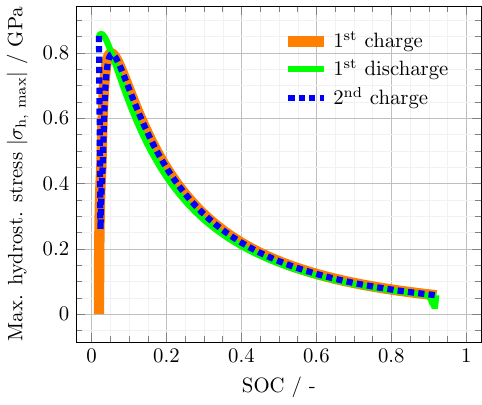}
    \caption{
      Absolute value of the maximal hydrostatic Cauchy
      stress~$\sigma_\text{h, max}$ over SOC for charging--discharging--charging
      cycle
      for the residual based error estimator
      with updated constants~$\gamma_1 = 0.001$ and $\gamma_2=1.0$.
    }
    \label{Schoof_Figure7:residual_stress_1d}
  \end{figure}
  \cref{Schoof_Figure7:residual_stress_1d}
  displays the development of the absolute value of the maximal hydrostatic
  stress~$\sigma_\text{h}$
  over the particle radius during the three half cycles.
  The hydrostatic Cauchy stress is defined
  as~$\sigma_\text{h} = (\sigma_\text{r} + 2 \sigma_{\phi})/3$
  with the stress in radial
  direction~$\sigma_\text{r}$
  and in tangential
  direction~$\sigma_{\phi}$.
  The evolution of the stress reveals high stresses at
  low~$\text{SOC}$
  values
  and lower stresses at high~$\text{SOC}$ values.
  The increase of the first charge at low~$\text{SOC}$
  results from the constant initial condition of the concentration.
  The small drop of the first discharge at high~$\text{SOC}$
  and the large drop of the second charge at low~$\text{SOC}$
  emerge from the change of stress curvature.
  For more information about this curvature change,
  see~\cite[Figure~7]{Schoof_Schoof_2023simulation}.
  In total,
  it can be concluded
  from the $\text{SOC}$ ranges and the behavior of the number of DOFs
  of~\cref{Schoof_Figure6:residual_dofs_1d}
  and~\cref{Schoof_Figure7:residual_stress_1d}
  that for the 1D spherical symmetric case
  larger maximal Cauchy stress values
  are related with
  a larger number of DOFs.

  \textbf{Example 2: 2D Quarter Ellipse.}
  \begin{figure}[t]
    \centering
    \includegraphics[width = 0.8\textwidth,
    page=1]{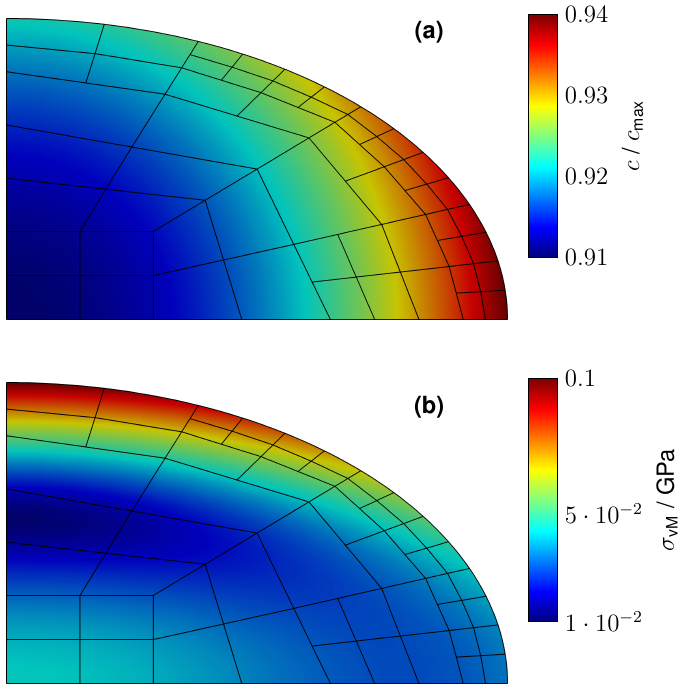}
    \caption{
      Numerical solution of the
      concentration~$c$
      in~(a)
      and
      von Mises
      stress~$\sigma_\text{vM}$
      in~(b)
      with the used mesh
      at the final simulation time~$t_\text{end}=0.9~\si{h}$
      for the updated version of the residual based error
      estimator
      with~$
      \gamma_1 = 0.001$
      and~$
      \gamma_2=1.0$.
    }
    \label{Schoof_Figure8:residual_2d}
  \end{figure}
  \begin{figure}[b]
    \centering
    \includegraphics[width = 0.77\textwidth,
    page=1]{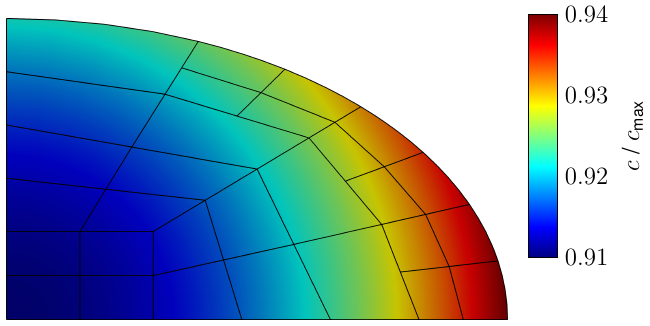}
    \caption{
      Numerical solution of the
      concentration~$c$
      with the used mesh
      at the final simulation time~$t_\text{end}=0.9~\si{h}$
      with the gradient recovery estimator.
    }
    \label{Schoof_Figure9:gradient_recovery_2d}
  \end{figure}
  The 2D computational domain is a proof of concept to show the feasibility of
  the gradient recovery estimator
  and updated version of the residual based error estimator
  as adaptive refinement strategies.
  For the simulations, we have modified some parameters in the following way:
  the
  number of initial refinements
  to~6,
  initial time step size
  to~$\num{1e-8}$
  and $t_\text{end} = 0.9$.
  In~\cref{Schoof_Figure8:residual_2d},
  the numerical solution of the
  concentration~$c$
  in Subfigure~(a)
  and of the von Mises
  stress~$\sigma_\text{vM} =\sqrt{\sigma_{11}^2+\sigma_{22}^2-
    \sigma_{11}\sigma_{22} + 3 \sigma_{12}^2}$
  in Subfigure~(b)
  are shown
  with the underlying mesh of the updated version of the residual based error
  estimator
  with~$\gamma_1 = 0.001$
  and~$\gamma_2=1.0$.
  Due to the elliptical domain,
  more lithium gets inside or out of the host particle
  at the area with the highest
  surface to volume ratio which is at the particle surface of the larger half
  axis~\cite[Section~2.2.3]{Schoof_Huttin_2014phase-field}.
  We point out
  that the mesh has higher refinements at areas of higher concentration values
  although the highest stresses are at the particle surface of the smaller half
  axis.
  The large von Mises stress~$\sigma_\text{vM}$
  at the smaller half axis occurs due to
  the forced zero stress in normal direction
  and the large tangential stresses in this area.
  This indicates
  that the concentration~$c$ has a larger influence than the von Mises
  stress~$\sigma_\text{vM}$
  on the spatial mesh strategy
  in the two-dimensional distribution at one time step.

  The numerical solution of the
  concentration~$c$
  with the underlying mesh
  of the gradient recovery refinement strategy is displayed
  in~\cref{Schoof_Figure9:gradient_recovery_2d}.
  Here,
  also a less refined mesh fulfills the tolerances
  of the adaptive solution algorithm
  compared to the updated version of the residual based strategy.
  This is similar to the findings of the 1D spherical symmetric case.



  \section{Conclusion and Outlook}
  \label{Schoof_sec:conclusion}

  \textbf{Conclusion.}
  In this work, three different spatial refinement strategies are investigated
  which have been applied to a chemo-elastic model for Li-diffusion within \aSi\
  particles. The ready to use Kelly error estimator~\cite{Schoof_Arndt_2023deal}
  of \texttt{deal.II}, the gradient recovery estimator
  of~\cite{Schoof_Castelli_2021efficient} and a residual based error estimator
  are
  compared. For the usage of the residual based error estimator, the strong
  formulation of the residual is derived explicitly in~\cref{Schoof_sec:theory}
  including many dependencies of our chemo-mechanical coupled model approach.
  In~\cref{Schoof_sec:results}, we considered the numerical results for a 3D
  spherical domain, reduced to 1D computational domain, and a 3D elliptical
  nanowire, reduced to 2D computational domain. For the 1D computational domain,
  we found out that the Kelly error estimator over-estimates the error at the
  boundary and therefore cannot be
  recommended for our application case. The
  gradient recovery estimator features higher refinement levels during the
  change
  of the lithium flux and rather lower refinement levels during solely the
  Li-diffusion periods. The residual based error estimator shows a stair-like
  behavior from the particle center to the particle surface. However, the cell
  error part leads to significant larger refinement levels compared to the face
  error part. Therefore, we have introduced an
  updated residual based error
  estimator with~$\gamma_1=0.001$
  (suitably chosen)
  and~$\gamma_2=1.0$, which reduces the influence
  of the cell error part. The updated version leads to a more efficient error
  estimator, since the number of DOFs is smaller by a factor of $10$ and the
  updated error estimator results in a comparable $L^2$- and~$H^1$-error
  compared to the residual case with~$\gamma_1=1.0$ and~$\gamma_2=1.0$.

  During a charging and discharging process,
  higher Cauchy stress values are related with a higher number of DOFs,
  whereas in the 2D computational domain, the distribution of Li-concentration
  has a larger influence on the mesh distribution
  compared to the Cauchy stress distribution at one specific time step.

  \textbf{Outlook.}
  Since we pointed out the significant influence of the cell error fraction of
  the residual based error estimator,
  further work is required for
  a closer examination and detailed
  calibration of the ratio of the cell and face
  errors of the residual based error estimator
  or for a potential influence of
  the finite element order.
  However, as charging and discharging is a highly dynamic process,
  a suitable calibration should be carefully checked.
  Another application of great interest is
  the consideration of another definition of the elastic
  strain
  tensor~$\te{E}_\el$
  as
  in~\cite{Schoof_Schoof_2023efficient}
  leading to an easy usage for further elasto-plastic couplings.
  In addition,
  the investigation of other battery active materials
  like graphite-silicon composites
  or
  materials for sodium-ion batteries
  is possible.
  Other 2D geometries and 3D geometries could also be taken into account.


  \section*{Declaration of competing interest}

  \noindent
  The authors declare that they have no known competing financial interests or
  personal relationships that could have appeared to influence the work in this
  paper.

  \section*{CrediT authorship contribution statement}
  \noindent
  \textbf{R. Schoof:}
  Conceptualization,
  Data curation,
  Formal Analysis,
  Investigation,
  Methodology,
  Software,
  Validation,
  Visualization,
  Writing -- original draft
  \noindent
  \textbf{L. Fl\"ur:}
  Formal Analysis,
  Methodology,
  Software,
  Visualization,
  Writing -- review \& editing
  \noindent
  \textbf{F. Tuschner:}
  Data curation,
  Writing -- review \& editing
  \noindent
  \textbf{W. D\"orfler:}
  Conceptualization,
  Funding acquisition,
  Project administration,
  Resources,
  Supervision,
  Writing -- review \& editing

  \section*{Acknowledgement}

  \noindent
  The authors thank
  G.~F.~Castelli for the software basis,
  L.~von Kolzenberg
  and L.~K\"obbing for intensive and constructive discussions about modeling
  silicon particles
  and T.~Laufer for proofreading.
  R.~S. and F.~T. acknowledge financial support by the German Research
  Foundation~(DFG) through the Research Training Group 2218
  SiMET~--~Simulation of Mechano-Electro-Thermal processes in Lithium-ion
  Batteries, project number 281041241.
  The authors acknowledge support by the state of Baden-Württemberg through
  bwHPC.

  \section*{Data Availability Statement}

  \noindent
  Data will be made available on request.

  \section*{Keywords}

  \noindent
  adaptive finite element method,
  residual based error estimator,
  finite deformation,
  lithium-ion batteries,
  numerical simulation


  \bibliography{Schoof_SchoofEtAl2024Residual}%


  \appendix

  \addcontentsline{toc}{section}{Appendices}
  \section*{Appendices}

  \renewcommand{\thesection}{A}
  \section{Abbreviations and Symbols}
  \label{app:abbreviations_and_symbols}

  \hspace*{-0.60cm}
  \vspace{-0.4cm}
  \begin{longtable}[l]{p{2.27cm}p{5.5cm}}
    \multicolumn{2}{l}{\textbf{Abbreviations}} \\
    \ \\
    AD & automatic differentiation \\
    \aSi & amorphous silicon \\
    \si{C}-rate & charging rate \\
    DAE & differential algebraic equation \\
    DOF & degree of freedom \\
    NDF & numerical differentiation formula \\
    OCV & open-circuit voltage \\
    pmv & partial molar volume \\
    SOC & state of charge
  \end{longtable}

  \hspace*{-0.60cm}
  \vspace{-0.4cm}
  \begin{longtable}[l]{p{2.27cm}p{10.5cm}}
    \textbf{Symbol} & \textbf{Description} \\
    \ \\
    \multicolumn{2}{l}{Latin symbols}
    \vspace{0.2cm}\\
    $c$ & concentration\\
    $\te{C}_\el$ & right Cauchy-Green tensor \\
    $\bm{\mathds{C}}$ & fourth-order stiffness tensor \\
    $E$ & Young's modulus\\
    $\te{E}_\el$ & elastic strain tensor \\
    $\te{F}
    $ & deformation gradient \\
    $\te{F} = \te{F}_\ch\te{F}_\el$ & multiplicative
    decomposition of $\te{F}$ \\
    $\te{F}_\ch$ & chemical deformation gradient \\
    $\te{F}_\el$ & elastic deformation gradient \\
    $\mathrm{Fa}$ & Faraday constant \\
    $G$ & shear modulus / second Lam\'{e} constant \\
    $\te{Id}$ & identity tensor \\
    $J = J_\ch J_\el$ & multiplicative
    decomposition of volume change \\
    $k_n$ & order of time integration algorithm at time $t_n$ \\
    ${m}$ & scalar valued mobility \\
    $\ve{n}$, $\ve{n}_0$ & normal vector on $\Omega$, $\Omega_0$ \\
    $\ve{N}$ & lithium flux \\
    ${N}_\text{ext}$ & external lithium flux \\
    $\te{P}$ & first Piola--Kirchhoff stress tensor \\
    $t$ & time \\
    $t_n$ & time step \\
    $t_\text{cycle}$ & cycle time \\
    $U_\text{OCV}$ & OCV curve \\
    $\ve{u} = \ve{x} - \ve{X}_0$ & time dependent displacement vector \\
    $v_\text{pmv}$ & partial molar volume of lithium \\
    $\ve{x}
    $ & position in Eulerian
    domain \\
    $\ve{X}_0$ & initial placement \\
    & \\
    %
    \multicolumn{2}{l}{Greek symbols} \\
    $\gamma_1$
    &
    factor for the cell error of the residual based error estimator \\
    $\gamma_2$
    &
    factor for the face error of the residual based error estimator \\
    $\Gamma_\text{in}$ & inner artificial boundary part \\
    $\lambda$ & first Lam\'{e} constant \\
    $\lambda_\ch$ & factor of concentration induced deformation gradient \\
    $\mu$ & chemical potential \\
    $\nu$ & Poisson's ratio \\
    $\Omega$ & Eulerian domain \\
    $\Omega_0$ & Lagrangian domain \\
    $\psi=\psi_\ch+\psi_\el$ & total energy \\
    $\psi_\ch$ & chemical part of total energy \\
    $\psi_\el$ & elastic part of total energy \\
    $\rho$ & density \\
    $\te{\sigma}$ & Cauchy stress tensor\\
    $\tau_n$ & time step size at time $t_n$\\
    & \\
    %
    \multicolumn{2}{l}{Mathematical symbols}	\\
    $\partial \Omega_0$ & boundary of $\Omega_0$ \\
    $\gradL$ & gradient vector in Lagrangian domain \\
    $\square \!:\! \tilde{\square}$ & reduction of two dimensions of two tensors
    $\square$ and $\tilde{\square}$\\
    $\bm{\mathcal{A}} [ \te{B}]$ &
    reduction of the last two dimensions of the third order
    tensor~$\bm{\mathcal{A}}$
    and the second order tensor $\te{B}$
    \\
    $\bm{\mathds{A}} [ \te{B}]$ &
    reduction of the last two dimensions of the fourth order
    tensor~$\bm{\mathds{A}}$
    and the second order tensor $\te{B}$
    \\
    $\partial_\square$ & partial derivative with respect to $\square$ \\
    & \\
    %
    %
    \multicolumn{2}{l}{Indices} \\
    ${\square}_{0}$ & considering variable in Lagrangian domain or initial
    condition
    \\
    ${\square}_\text{ch}$ & chemical part of~$\square$ \\
    ${\square}_{\text{com}}$ & computational part of $\square$\\
    ${\square}_\text{el}$ & elastic part of~$\square$ \\
    ${\square}_{\text{max}}$ & maximal part of $\square$\\
    ${\square}_{\text{min}}$ & minimal part of $\square$\\
  \end{longtable}

  \renewcommand{\thesection}{B}
  \section{Simulation Parameters}
  \label{Schoof_app:simulation_parameters}
  The dimensionless variables of the considered model equations are given
  in~\cref{Schoof_Table2:normalization}
  and
  the used model parameters for the numerical simulations are listed
  in~\cref{Schoof_Table3:parameters}.

  \setcounter{table}{1}

  \begin{table}[b]
    \centering
    \caption{Dimensionless variables of the used model equations.}
    \label{Schoof_Table2:normalization}
    \begin{tabular}[t]{@{}lllll@{}}
      \toprule
      $\tilde{\ve{X}}_0 = \ve{X}_0 / L_0$ \quad \quad
      & $\tilde{t} = t / t_\text{cycle}$ \quad \quad
      & $\tilde{\ve{u}} = \ve{u} / L_0$ \quad \quad \quad
      & $\tilde{c} = c / c_{\max}$ \quad \quad
      & $\tilde{v}_\text{pmv} = {v}_\text{pmv} c_{\max}$
      \\
      \multicolumn{2}{l}{
        $\widetilde{\rho \psi} = \rho \psi / R_\text{gas} T c_{\max}$}
      & $\tilde{\mu} = \mu / R_\text{gas} T$
      & \multicolumn{2}{l}{$\tilde{E} = E / R_\text{gas} T c_{\max}$}
      \\
      \multicolumn{3}{l}{
        $\mathrm{Fo} = D t_\text{cycle}/ L_0^2$
        \quad \quad
        $\tilde{N}_{\ext} = {N}_{\ext} t_\text{cycle}/ L_0 c_{\max}$}
      \quad
      &
      \multicolumn{2}{l}{
        $\tilde{U}_\text{OCV}
        = \mathrm{Fa} \, {U}_\text{OCV} / R_\text{gas} T$}
      \\
      \bottomrule
    \end{tabular}
  \end{table}

  \begin{table}[h]
    \centering
    \caption{Model parameters for numerical experiments
      \cite{Schoof_Kolzenberg_2022chemo-mechanical}.
    }
    \label{Schoof_Table3:parameters}
    \begin{tabular}[t]{@{}llccc@{}}
      \toprule
      \textbf{Description} & \textbf{Symbol} &
      \multicolumn{1}{c}{\textbf{Value}} & \textbf{Unit} &
      \multicolumn{1}{c}{\textbf{Dimensionless}} \\

      \midrule

      Universal gas constant & $R_\text{gas}$ & $8.314$ &
      \si{\joule\per\mol\per\kelvin}
      & $1$ \\[\defaultaddspace]

      Faraday constant & $\mathrm{Fa}$ & $96485$ & \si{\joule\per\volt\per\mol}
      &
      $1$
      \\[\defaultaddspace]

      Operation temperature & $T$ & $298.15$ & \si{\kelvin} & $1$
      \\[\defaultaddspace]

      \midrule

      \multicolumn{5}{c}{Silicon} \\

      \midrule

      Particle length scale & $L_0$ & $\num{50e-9}$ &
      \si{\meter} & 1
      \\[\defaultaddspace]

      Diffusion coefficient & $D$ & $\num{1e-17}$ &
      \si{\square\meter\per\second} &
      $14.4$ \\[\defaultaddspace]

      OCV curve & $U_\text{OCV}$ &
      \cite{Schoof_Kolzenberg_2022chemo-mechanical} &
      \si{\volt}
      &
      \cite{Schoof_Kolzenberg_2022chemo-mechanical}
      \\[\defaultaddspace]

      Young's modulus & $E$ & $\num{90.13e9}$ & \si{\pascal} & $116.74$
      \\[\defaultaddspace]

      Partial molar volume & $v_\text{pmv}$ & $\num{10.96e-6}$ &
      \si{\cubic\meter\per\mol} & $3.41$
      \\[\defaultaddspace]

      Maximal concentration & $c_\text{max}$ & $\num{311.47e3}$ &
      \si{\mol\per\cubic\meter} & $1$
      \\[\defaultaddspace]

      Initial concentration & $c_0$ & $\num{6.23e3}$ &
      \si{\mol\per\cubic\meter} & $\num{2e-2}$
      \\[\defaultaddspace]

      Poisson's ratio & $\nu$ & $0.22$ & \si{-} & $0.22$
      \\[\defaultaddspace]

      \bottomrule
    \end{tabular}
  \end{table}


\end{document}